\documentclass{amsart}
\usepackage{hyperref}

\newtheorem{theorem}{Theorem}[section]
\newtheorem{definition}[theorem]{Definition}
\newtheorem{lemma}[theorem]{Lemma}
\newtheorem{corollary}[theorem]{Corollary}

\newtheorem{hypothesis}[theorem]{Hypothesis {\bf H.}\hspace*{-0.6ex}}

\newcommand{\R}{{\mathbb R}}

\newcommand{\C}{{\mathbb C}}

\newcommand{\nn}{\nonumber}
\newcommand{\be}{\begin{equation}}
\newcommand{\ee}{\end{equation}}
\newcommand{\bea}{\begin{eqnarray}}
\newcommand{\eea}{\end{eqnarray}}
\newcommand{\ul}{\underline}
\newcommand{\ol}{\overline}
\newcommand{\ti}{\tilde}

\newcommand{\spr}[2]{\langle #1 , #2 \rangle}

\newcommand{\sgn}{\mathrm{sgn}}
\newcommand{\tr}{\mathrm{tr}}

\newcommand{\dom}{\mathfrak{D}}
\newcommand{\fdom}{\mathfrak{Q}}

\DeclareMathOperator{\Ran}{Ran}

\newcommand{\floor}[1]{\lfloor#1 \rfloor}
\newcommand{\ceil}[1]{\lceil#1 \rceil}
\newcommand{\norm}[1]{\lVert#1 \rVert}
\newcommand{\abs}[1]{\lvert#1 \rvert}

\newcommand{\eps}{\varepsilon}

\newcommand{\sig}{\sigma}
\newcommand{\lam}{\lambda}

\numberwithin{equation}{section}

\begin{document}

\title[Relative Oscillation Theory Extended]{Relative Oscillation Theory for Sturm--Liouville Operators Extended}

\author[H. Kr\"uger]{Helge Kr\"uger}
\address{Faculty of Mathematics\\
Nordbergstrasse 15\\ 1090 Wien\\ Austria}
\email{\href{mailto:helge.krueger@rice.edu}{helge.krueger@rice.edu}}
\urladdr{\href{http://www.mat.univie.ac.at/~helge/}{http://www.mat.univie.ac.at/\~{}helge/}}
\curraddr{Department of Mathematics, Rice University, Houston, TX 77005, USA}

\author[G. Teschl]{Gerald Teschl}
\address{Faculty of Mathematics\\
Nordbergstrasse 15\\ 1090 Wien\\ Austria\\ and International Erwin Schr\"odinger
Institute for Mathematical Physics, Boltzmanngasse 9\\ 1090 Wien\\ Austria}
\email{\href{mailto:Gerald.Teschl@univie.ac.at}{Gerald.Teschl@univie.ac.at}}
\urladdr{\href{http://www.mat.univie.ac.at/~gerald/}{http://www.mat.univie.ac.at/\~{}gerald/}}

\thanks{{\it To appear in J. Funct. Anal.}}
\thanks{{\it Research supported by the Austrian Science Fund (FWF) under Grant No.\ Y330}}


\begin{abstract}
We extend relative oscillation theory to the case of Sturm--Liouville operators
$H u = r^{-1}(-(pu')'+q u)$ with different $p$'s. We show that the weighted number of zeros
of Wronskians of certain solutions equals the value of Krein's spectral shift function inside
essential spectral gaps.
\end{abstract}

\maketitle

\section{Introduction}

In \cite{kt} we have developed an analog of classical oscillation theory for Sturm--Liouville
operators which, rather than measuring the spectrum of one single 
operator, measures the difference between the spectra of two different operators.
Hence the name relative oscillation theory. The main idea behind this extension is to
replace zeros of solutions of one operator by weighted zeros of Wronskians of solutions of
two different operators. That zeros of the Wronskian are related to oscillation theory is indicated
by an old paper of Leighton \cite{lei}, who noted that if two solutions have a
non-vanishing Wronskian, then their zeros must intertwine each other. Their
use as an adequate tool for the investigation of the spectrum of one single operator is due to
Gesztesy, Simon, and one of us \cite{gst}.

The purpose of this paper is to extend relative oscillation theory
for two different Sturm--Liouville equations
\be
\tau_j = \frac{1}{r} \Big(- \frac{d}{dx} p_j \frac{d}{dx} + q_j \Big), \qquad j=0,1.
\ee
In \cite{kt} we considered the case $p_0=p_1$, here we want to extend relative
oscillation theory to the case $p_0\ne p_1$. In particular, for $H_j$, $j=0,1$, self-adjoint
operators associated with $\tau_j$, we want to show that the
weighted number of zeros of Wronskians of certain solutions equals the value of Krein's
spectral shift function $\xi(\lam,H_1,H_0)$ inside essential spectral gaps. To do this, and
to make sure that the spectral shift function is well-defined, we will need to find a continuous
path connecting the operators $H_0$ and $H_1$ in the metric introduced by the trace norm
of resolvent differences.

In Section~\ref{sec:wz} we will recall the necessary background and fix our notation.
Moreover, we will present the basic result for the case of regular operators.
In Section~\ref{sec:sturthm} we have a quick look at Sturm's classical comparison theorem
for zeros of solutions and its extension to zeros of Wronskians of solutions.
Section~\ref{sec:rot} is concerned with relative oscillation theory for singular operators
and contains our key result, Theorem~\ref{thmsing}, which connects the weighted zeros of
Wronskians with Krein's spectral shift function. The remaining sections contain the proofs
for our main results and our final appendix collects some facts on the spectral shift functions
plus some abstract results which form the functional analytic core of the proof of our
main theorem.

\section{Weighted zeros of Wronskians, Pr\"ufer angles, and regular operators}
\label{sec:wz}

We begin by fixing our notation and reviewing some simple facts from \cite{kt}.
In particular, we refer to \cite{kt} for further details.

We will consider Sturm--Liouville operators on $L^{2}((a,b), r\,dx)$
with $-\infty \le a<b \le \infty$ of the form
\begin{equation} \label{stli}
\tau = \frac{1}{r} \Big(- \frac{d}{dx} p \frac{d}{dx} + q \Big),
\end{equation}
where the coefficients $p,q,r$ are real-valued
satisfying
\begin{equation}
p^{-1},q,r \in L^1_{loc}(a,b), \quad p,r>0.
\end{equation}
We will use $\tau$ to describe the formal differentiation expression and
$H$ for the operator given by $\tau$ with separated boundary conditions at
$a$ and/or $b$.

If $a$ (resp.\ $b$) is finite and $q,p^{-1},r$ are in addition integrable
near $a$ (resp.\ $b$), we will say $a$ (resp.\ $b$) is a \textit{regular}
endpoint.  We will say $\tau$ respectively $H$ is \textit{regular} if
both $a$ and $b$ are regular. 

For every $z\in\C\backslash\sig_{ess}(H)$ there is a unique (up to a constant) solution
$\psi_-(z,x)$ of $\tau u = z u$ which is in $L^2$ near $a$ and satisfies the
boundary condition at $a$ (if any). Similarly there is such a solution $\psi_+(z,x)$
near $b$.

One of our main objects will be the (modified) Wronskian
\be
W_x(u_0,u_1)= u_0(x)\, p_1(x)u_1'(x) - p_0(x)u_0'(x)\, u_1(x)
\ee
of two functions $u_0$, $u_1$ and its zeros. Here we think of $u_0$ and $u_1$ as two
solutions of two different Sturm--Liouville equations
\be \label{stlij}
\tau_j = \frac{1}{r} \Big(- \frac{d}{dx} p_j \frac{d}{dx} + q_j \Big), \qquad j=0,1.
\ee
Under these assumptions $W_x(u_0,u_1)$ is absolutely continuous and satisfies
\be\label{dwr}
W'(u_0, u_1) = (q_1 - q_0) u_0 u_1 +
\left(\frac{1}{p_0} - \frac{1}{p_1}\right) p_0 u_0' p_1 u_1'.
\ee
Next we recall the definition of Pr\"ufer variables $\rho_u$, $\theta_u$
of an absolutely continuous function $u$:
\be\label{eq:defprue}
u(x)=\rho_u(x)\sin(\theta_u(x)), \qquad
p(x) u'(x)=\rho_u(x) \cos(\theta_u(x)).
\ee
If $(u(x), p(x) u'(x))$ is never $(0,0)$ and $u, pu'$ are absolutely continuous, then $\rho_u$ is
positive and $\theta_u$ is uniquely determined once a value of
$\theta_u(x_0)$ is chosen by requiring continuity of $\theta_u$.

Notice that
\be \label{wpruefer}
W_x(u,v)= -\rho_u(x)\rho_v(x)\sin(\Delta_{v,u}(x)), \qquad
\Delta_{v,u}(x) = \theta_v(x)-\theta_u(x).
\ee
Hence the Wronskian vanishes if and only if the two Pr\"ufer angles differ by
a multiple of $\pi$. We will call the total difference
\be
\#_{(c,d)}(u_0,u_1) = \ceil{\Delta_{1,0}(d) / \pi} - \floor{\Delta_{1,0}(c) / \pi} -1
\ee
the number of weighted sign flips in $(c,d)$, where we have written $\Delta_{1,0}(x)=
\Delta_{u_1,u_0}$ for brevity.

We take two solutions $u_j$, $j=1,2$, of $\tau_j u_j =\lam_j u_j$ and
associated Pr\"ufer variables $\rho_j$, $\theta_j$. Since we can
replace $q \to q - \lam r$ it is no restriction to assume $\lam_0=\lam_1=0$.
We remark, that in (\ref{eq:defprue}) one has to take $p_j$ as $p$ for $u_j$,
$j=0,1$.

\begin{lemma}\label{lem:delinc}
Abbreviate $\Delta_{1,0}(x) = \theta_1(x) - \theta_0(x)$ and
suppose $\Delta_{1,0}(x_0) \equiv 0 \mod \pi$.
If $q_0(x)-q_1(x)$ and $p_0(x)-p_1(x)$ are
(i) negative, (ii) zero, or (iii) positive for a.e.
$x \in (x_0, x_0 + \eps)$ respectively for a.e.
$x \in (x_0 - \eps, x_0)$, then the same holds true
for $(\Delta_{1,0}(x) - \Delta_{1,0}(x_0))/(x - x_0)$.
\end{lemma}

\begin{proof}
By (\ref{dwr}) we have
\begin{align} \nn
W_x(u_0,u_1) &= - \rho_0(x)\rho_1(x)\sin(\Delta_{1,0}(x))\\
&=  - \int_{x_0}^x \Big( (q_0(t)-q_1(t)) u_0(t) u_1(t) +
(\frac{1}{p_1(t)}-\frac{1}{p_0(t)}) p_0 u_0'(t) p_1 u_1'(t) \Big)dt.
\end{align}
Case (ii) follows. For (i) and (iii), first note that
if $u_j(x_0) = 0$, $j = 0,1$, we have that
$u_j$ and $p_j u_j'$, $j=0,1$ have the same sign close to $x_0$,
and thus the result follows.

Now, look at $P(u_0,u_1)=\frac{u_0}{u_1} W(u_0,u_1)$ (compare (\ref{eq:pic}) below)
(resp.\ $P(u_1,u_0)$) and note that $u_0/u_1$ has constant sign near $x_0$.
The result now follows using the fact that the
derivate $P'(u_0,u_1)$ is always negative by the Picone identity (\ref{eq:pd}) below.
\end{proof}

\noindent
Hence $\#_{(c,d)}(u_0,u_1)$ counts the weighted sign flips of the Wronskian
$W_x(u_0,u_1)$, where a sign flip is counted as $+1$ if $q_0-q_1$ and $p_0-p_1$ are positive
in a neighborhood of the sign flip, it is counted as $-1$ if $q_0-q_1$ and $p_0-p_1$ are negative
in a neighborhood of the sign flip. In particular, we obtain

\begin{lemma}\label{lemnbzer2}
Let $u_0$, $u_1$ solve $\tau_j u_j = 0$, $j=0,1$, where $p_0-p_1\ge 0$ and $q_0-q_1\ge 0$.
Then $\#_{(a,b)}(u_0,u_1)$ equals the number sign flips of $W(u_0,u_1)$
inside the interval $(a,b)$.
\end{lemma}

\noindent
Finally, we have the following extension of \cite[Thm.~2.3]{kt} to the case $p_0\ne p_1$.

\begin{theorem} \label{thm:reg}
Let $H_0$, $H_1$ be regular Sturm--Liouville operators associated with (\ref{stlij}) and
the same boundary conditions at $a$ and $b$. Then
\be \label{eqreg}
\dim\Ran\, P_{(-\infty, \lam_1)}(H_1) - \dim\Ran\, P_{(-\infty, \lam_0]}(H_0) =
\#_{(a,b)}(\psi_{0,\pm}(\lam_0), \psi_{1,\mp}(\lam_1)).
\ee
\end{theorem}

\noindent
The proof will be given in Section~\ref{sec:proofreg}.

\section{Sturm's comparison theorem}
\label{sec:sturthm}

One of the core ingredients of oscillation theory is Sturm's comparison theorem
for zeros of solutions. We begin by recalling this classical result.

Let $u_j$ solve $\tau_j u_j = \lam_j u_j$, where without loss of generality
we assume $\lam_0=\lam_1=0$. For $x$ with $u_1(x)\neq 0$ we introduce
\be\label{eq:pic}
P_x(u_0,u_1) = \frac{u_0(x)}{u_1(x)} W_x(u_0, u_1) =
-\rho_0^2(x) \frac{\sin(\theta_0(x))\sin(\Delta_{1,0}(x))}{\sin(\theta_1(x))}.
\ee
Obviously $P(u_0,u_1)$ is zero if either $u_0$ or the Wronskian $W(u_0, u_1)$
vanishes. Moreover, a straightforward computation, verifies the Picone identity
(see \cite[(2.6.4)]{zet})
\be\label{eq:pd}
P'(u_0,u_1) = (q_1-q_0) u_0^2 + (p_1-p_0)u_0'^2
-p_1\left(u_0' -  \frac{u_0 u_1'}{u_1}\right)^2,
\ee
which shows that $P(u_0,u_1)$
is a nonincreasing function if $q_1 \leq q_0$ and $0 < p_1 \leq p_0$.

\begin{theorem}[Sturm's Comparison theorem]\label{thm:sturclaform}
Let $q_0-q_1 \geq 0$, $p_0 -p_1\geq 0$, with once
strict inequality, and $\tau_j u_j = 0$, $j=0,1$.
Then between any two zeros of $u_0$ or $W(u_0, u_1)$,
there is a zero of $u_1$.

Similarly, between two zeros of $u_1$, which are not at the same time zeros of $u_0$,
there is at least one zero of $u_0$ or $W(u_0,u_1)$.
\end{theorem}

\begin{proof}
Assume that $u_1$ has no zero,
$P(u_0,u_1)$ would be well defined
on the closed interval between the zeros, and be zero
at its end points. This contradicts monotonicity of $P(u_0,u_1)$.
The second claim is similar.
\end{proof}

\noindent
Note that this version is slightly more general then the one usually found in the
literature (cf., e.g, \cite{zet}) since it includes the case of zeros of Wronskians. For the
case $p_0 = p_1$ this was already pointed out in \cite{gst}. Moreover, in this case one
can also allow zeros of the Wronskian at singular endpoints \cite[Cor.~2.3]{gst}.

Next, the comparison theorem for Wronskians from \cite{kt} carries over to the case
$p_0\ne p_1$ without modifying the proof.

\begin{theorem}[Comparison theorem for Wronskians]\label{thmscw}
Suppose $u_j$ satisfies $\tau_j u_j = \lam_j u_j$, $j=0,1,2$,
where $\lam_0 r -q_0 \le \lam_1 r - q_1 \le \lam_2 r - q_2$,
$p_0 \geq p_1 \geq p_2$.

If $c<d$ are two zeros of $W_x(u_0,u_1)$ such that $W_x(u_0,u_1)$ does not
vanish identically, then there is at least one sign flip of $W_x(u_0,u_2)$ in $(c,d)$.
Similarly, if $c<d$ are two zeros of $W_x(u_1,u_2)$ such that $W_x(u_1,u_2)$ does not
vanish identically, then there is at least one sign flip of $W_x(u_0,u_2)$ in $(c,d)$.
\end{theorem}

\section{Relative Oscillation Theory}
\label{sec:rot}

After these preparations we are now ready to extend relative oscillation theory 
to the case $p_0\ne p_1$. Except for Lemma~\ref{lem:nonoscingap2} and
our key result Theorem~\ref{thmsing}, all results in this section are straightforward
modifications of the analog results in \cite{kt} and hence we omit
the corresponding proofs.

\begin{definition}\label{def:wsf}
For $\tau_0$, $\tau_1$ possibly singular Sturm--Liouville operators as in (\ref{stlij}) on $(a,b)$,
we define
\be
\underline{\#}(u_0,u_1) = \liminf_{d \uparrow b,\,c \downarrow a} \#_{(c,d)}(u_0,u_1)
\quad\mbox{and}\quad
\overline{\#}(u_0,u_1) = \limsup_{d \uparrow b,\,c \downarrow a} \#_{(c,d)}(u_0,u_1),
\ee
where $\tau_j u_j = \lam_j u_j$, $j=0,1$.

We say that $\#(u_0,u_1)$ exists, if
$\overline{\#}(u_0,u_1)=\underline{\#}(u_0,u_1)$, and write
\be
\#(u_0,u_1) = \overline{\#}(u_0,u_1)=\underline{\#}(u_0,u_1).
\ee
in this case.
\end{definition}

\noindent
By Lemma~\ref{lem:delinc} one infers that $\#(u_0,u_1)$ exists
if $p_0-p_1$ and $q_0-\lam_0r- q_1+\lam_1r$ have the same definite sign
near the endpoints $a$ and $b$.

\begin{theorem}[Triangle inequality for Wronskians]\label{thm:wtriang}
Suppose $u_j$, $j=0,1,2$ are given functions with $u_j$, $p_ju_j'$ absolutely continuous and
$(u_j(x), p_j(x) u_j'(x)) \neq (0,0)$ for all $x$. Then
\be 
\ul{\#}(u_0,u_1) + \ul{\#}(u_1,u_2) - 1 \leq \ul{\#}(u_0,u_2) \leq
\ul{\#}(u_0,u_1) + \ul{\#}(u_1,u_2) + 1
\ee
and similarly for $\ul{\#}$ replaced by $\ol{\#}$.
\end{theorem}

\noindent
We recall that in classical oscillation theory
$\tau$ is called oscillatory if a solution of $\tau u = 0$ has infinitely many
zeros.

\begin{definition}\label{def:relosc}
We call $\tau_1$ relatively nonoscillatory with respect to
$\tau_0$, if the quantities $\underline{\#}(u_0, u_1)$ and
$\overline{\#}(u_0, u_1)$ are finite for all solutions
$\tau_j u_j = 0$, $j = 0,1$.
We call $\tau_1$ relatively oscillatory with respect to
$\tau_0$, if one of the quantities $\underline{\#}(u_0, u_1)$ or
$\overline{\#}(u_0, u_1)$ is infinite for some solutions
$\tau_j u_j = 0$, $j = 0,1$.
\end{definition}

\noindent
Note that this definition is in fact independent of the solutions
chosen as a straightforward application of our triangle inequality
(cf.\ Theorem~\ref{thm:wtriang}) shows.

\begin{corollary}\label{cor:indisol}
Let $\tau_j u_j = \tau_j v_j = 0$, $j=0,1$. Then
\be
|\underline{\#}(u_0, u_1) - \underline{\#}(v_0, v_1)| \le 4,\quad
|\overline{\#}(u_0, u_1) - \overline{\#}(v_0, v_1)| \le 4.
\ee
\end{corollary}

\noindent
The bounds can be improved using our comparison theorem for
Wronskians to be $\leq 2$ in the case of perturbations of definite sign.

If $\tau_0$ is nonoscillatory our definition reduces to the classical one.

\begin{lemma}
Suppose $\tau_0$ is a nonoscillatory operator, then
$\tau_1$ is relatively nonoscillatory (resp. oscillatory) with
respect to $\tau_0$, if and only if $\tau_1$ is
nonoscillatory (resp. oscillatory).
\end{lemma}

\noindent
To demonstrate the usefulness of Definition~\ref{def:relosc},
we now establish its connection with the spectra
of self-adjoint operators associated with $\tau_j$, $j=0,1$.

\begin{theorem} \label{thm:rosc}
Let $H_j$ be
self-adjoint operators associated with $\tau_j$, $j=0,1$. Then
\begin{enumerate}
\item
$\tau_0-\lam_0$ is relatively nonoscillatory with respect to $\tau_0-\lam_1$
if and only if $\dim\Ran P_{(\lam_0,\lam_1)}(H_0)<\infty$.
\item
Suppose $\dim\Ran P_{(\lam_0,\lam_1)}(H_0)<\infty$ and
$\tau_1-\lam$ is relatively nonoscillatory with respect to $\tau_0-\lam$ for one
$\lam \in [\lam_0,\lam_1]$.  Then it is relatively nonoscillatory for
all $\lam\in[\lam_0,\lam_1]$ if and only if $\dim\Ran P_{(\lam_0,\lam_1)}(H_1)<\infty$.
\end{enumerate}
\end{theorem}

\noindent
For a practical application of this theorem one needs of course criteria
when  $\tau_1-\lam$ is relatively nonoscillatory with respect to $\tau_0-\lam$
for $\lam$ inside an essential spectral gap. Without loss of generality we only
consider the case where one endpoint is regular.

\begin{lemma}\label{lem:nonoscingap2}
Let $H_0$ be bounded from below. Suppose $a$ is regular ($b$ singular) and
\begin{enumerate}
\item
$\lim_{x\to b} r(x)^{-1} (q_0(x) - q_1(x)) = 0$, $\frac{q_0}{r}$ is bounded near $b$, and
\item
$\lim_{x\to b} p_1(x) p_0(x)^{-1} = 1$.
\end{enumerate}
Then $\sig_{ess}(H_0)=\sig_{ess}(H_1)$ and
$\tau_1 - \lam$ is relatively nonoscillatory with respect
to $\tau_0 - \lam$ for every $\lam \in \R \backslash\sigma_{ess}(H_0)$.

The analogous result holds for $a$ singular and $b$ regular.
\end{lemma}

\noindent
The proof will be given in Section~\ref{sec:proofreg}. 

Our next task is to reveal the precise relation between the number of
weighted sign flips and the spectra of $H_1$ and $H_0$. The special
case $H_0=H_1$ is covered by \cite{gst}:

\begin{theorem}[\cite{gst}] \label{thm:gst}
Let $H_0$ be a self-adjoint operator associated with $\tau_0$ and suppose
$[\lam_0,\lam_1]\cap\sig_{ess}(H_0)=\emptyset$.
Then
\be
 \dim\Ran P_{(\lam_0, \lam_1)} (H_0) = \#(\psi_{0,\mp}(\lam_0),\psi_{0,\pm}(\lam_1)).
\ee
\end{theorem}

\noindent
Combining this result with our triangle inequality already gives some rough
estimates.

\begin{lemma}\label{lem:estproj}
Let $H_0$, $H_1$ be self-adjoint operators associated with $\tau_0$, $\tau_1$,
respectively, and separated boundary conditions. Suppose 
that $(\lam_0, \lam_1)\subseteq\R\backslash(\sig_{ess}(H_0)\cup\sig_{ess}(H_1))$, then
\begin{align}
\nonumber
\dim\Ran P_{(\lam_0, \lam_1)} (H_1) &-
\dim\Ran P_{(\lam_0, \lam_1)} (H_0)\\
& \leq \underline{\#}(\psi_{1,\mp} (\lam_1), \psi_{0,\pm} (\lam_1)) - 
\overline{\#}(\psi_{1,\mp} (\lam_0), \psi_{0,\pm} (\lam_0)) + 2,
\end{align}
respectively,
\begin{align}
\nonumber
\dim\Ran P_{(\lam_0, \lam_1)} (H_1) &-
\dim\Ran P_{(\lam_0, \lam_1)} (H_0)\\
& \geq \overline{\#}(\psi_{1,\mp} (\lam_1), \psi_{0,\pm} (\lam_1)) - 
\underline{\#}(\psi_{1,\mp} (\lam_0), \psi_{0,\pm} (\lam_0)) - 2.
\end{align}
\end{lemma}

\noindent
To extend Theorem~\ref{thm:reg} to the singular case, we need to require
the following hypothesis similar to \cite[H.3.15]{kt}.

\begin{hypothesis} \label{hyp:h0h1}
Suppose $H_0$ and $H_1$ are self-adjoint operators
associated with $\tau_0$ and $\tau_1$
and separated boundary conditions (if any).
Introduce
\begin{align}\nn
A_0 &= \frac{1}{r} (r p_0)^{1/2} \frac{d}{dx}, \\ \label{defA0}
\dom(A_0) &= \{ f\in L^2((a,b),r\, dx)| f\in AC_{loc}(a,b), \: \sqrt{p_0}f'\in L^2(a,b)\}
\end{align}
\begin{enumerate}
\item $r^{-1} q_0$ is infinitesimally form bounded with respect to $A_0^* A_0$.
\item $r^{-1}(q_1 - q_0)$ is infinitesimally form bounded with respect to $H_0$.
\item There is a $C_1>1$ such that $C_1^{-1} \leq p_0(x)^{-1} p_1(x) \leq C_1$ for all $x$.
\item $r^{-1}\abs{r(p_0-p_1)}^{1/2} \frac{d}{dx} R_{H_0} (z)$
and $\abs{(r^{-1}(q_1 - q_0)}^{1/2} R_{H_0}(z)$
are Hilbert--Schmidt for one (and hence for all)
$z\in\rho(H_0)$.
\end{enumerate}
\end{hypothesis}

\noindent
We note that the conditions of the last hypothesis are for example
satisfied for periodic operators if the coefficients are continuous and
$p_0^{-1} - p_1^{-1}$ and $q_0-q_1$ are integrable.

It will be shown in Section~\ref{sec:xi} that these conditions ensure that we can interpolate
between $H_0$ and $H_1$ using operators $H_\eps$, $\eps\in[0,1]$, such that the resolvent difference
of $H_0$ and $H_\eps$ is continuous in $\eps$ with respect to the trace norm. Hence we can fix
the spectral shift function $\xi(\lam,H_1,H_0)$ by requiring $\eps \mapsto \xi(\lam,H_\eps,H_0)$
to be continuous in $L^1(\R,(\lam^2+1)^{-1}d\lam)$, where we of course set $\xi(\lam,H_0,H_0)=0$
(see Lemma~\ref{lem:resolvconv}). While $\xi$ is only defined a.e., it is constant on the intersection
of the resolvent sets $\R\cap\rho(H_0)\cap\rho(H_1)$, and we will require it to be continuous there.
In particular, note that by Weyl's theorem the essential spectra of $H_0$ and
$H_1$ are equal, $\sig_{ess}(H_0)=\sig_{ess}(H_1)$.

\begin{theorem} \label{thmsing}
Let $H_0$, $H_1$ satisfy Hypothesis~\ref{hyp:h0h1}. Then
for every $\lam\in \rho(H_0)\cap\rho(H_1)\cap\R$, we have
\be
\xi(\lam,H_1,H_0) =
\#(\psi_{0,\pm}(\lam), \psi_{1,\mp}(\lam)).
\ee
\end{theorem}

\section{Proofs of Lemma~\ref{lem:nonoscingap2} and the regular case}
\label{sec:proofreg}

To prove Lemma~\ref{lem:nonoscingap2}, we need the
following  modification of \cite[Lem.~3.9]{kt}:

\begin{lemma}
Let $(\lam_0,\lam_1)\subseteq\R\backslash\sigma_{ess}(H_0)$,
$\lam\in(\lam_0,\lam_1)$. If $p_0=p_1$ and $\lam_0< r^{-1}(q_1-q_0)-\lam<\lam_1$
(at least near singular endpoints), then $\tau_1-\lam$ is relatively nonoscillatory
with respect to $\tau_0-\lam$.
\end{lemma}

\begin{proof}
Using our comparison theorem, we have that from
$\#(u_0(\lam_0),u_0(\lam_1))<\infty$, we obtain
$$
\ol{\#}(u_0(\lam_1),u_1(\lam))<\infty,\quad
\ul{\#}(u_0(\lam_0),u_1(\lam))>-\infty
$$
now the result follows as in \cite[Lem.~3.9]{kt} by
$$
\ol{\#}(u_0(\lam),u_1(\lam))\leq\#(u_0(\lam),u_0(\lam_1))+\ol{\#}(u_0(\lam_1),u_1(\lam))+1
$$
as follows from the triangle inequality for Wronskians (\cite[Thm.~3.4]{kt}) and \cite[Thm.~3.8~(i)]{kt}.
\end{proof}

\noindent
Our next proof will require the following resolvent relation for form perturbations. It is a special
case from \cite[Sect.~VI.3]{kato} (see in particular equation (VI.3.10); compare also
Sect.~II.3. in \cite{shqf}).

\begin{lemma}\label{lem:secresoform}
Let $H_0$ be a self-adjoint operator which is bounded from below and let
$\lam$ be below its spectrum. Let $V$ be relatively form bounded with respect to $H_0$
and with bound less than one. Then, we have that
$H = H_0 + V$ is self-adjoint and for its resolvent we have
\be
R_H(z) = R_{H_0}^{1/2}(\lam) (1 - (z-\lam) R_{H_0}(\lam) + C)^{-1} R_{H_0}^{1/2}(\lam).
\ee
Here $C$ is the bounded operator associated with the quadratic form
\be
\psi\mapsto\spr{R_{H_0}^{1/2}(\lam) \psi}{ V R_{H_0}^{1/2}(\lam) \psi}.
\ee
\end{lemma}

\noindent
We remark, that here and in what follows sums of operators have to be understood as forms sums.
Now we come to the

\begin{proof}[Proof of Lemma~\ref{lem:nonoscingap2}]
We first show that $\sig_{ess}(H_0) = \sig_{ess}(H_1)$.
First of all, note that imposing an additional Dirichlet boundary condition at
some point $b_n\in(a,b)$ implies that the resolvents of the original and the
perturbed operator differ by a rank one perturbation (cf., e.g, \cite{wdln}). Furthermore, the
perturbed operator decomposes into a direct sum of two operators,
one regular part on $(a,b_n)$ and one singular part on $(b_n,b)$. Since the
resolvent of a regular Sturm--Liouville operator is Hilbert--Schmidt, the
only interesting part for the essential spectrum is the singular operator
on $(b_n,b)$. Denote the corresponding operators by $H_j^n$, $j=1,2$.
(i.e., $H_j^n$ is $H_j$ restricted to $(b_n,b)$ with a Dirichlet boundary
condition at $b_n$). Then it suffices to show that the resolvent difference of
$H_1^n$ and $H_0^n$ can be made arbitrarily small by choosing
$b_n$ close to $b$.

Recall the definition of $A_0$ from (\ref{defA0}) and
note that since $\frac{q_0}{r}$ is bounded (for $b_n$ sufficiently large),
$A_0 R_{H_0}^{1/2}(-\lam)$ is bounded for $-\lam<\sig(H_0)$.
By virtue of Lemma~\ref{lem:secresoform} we conclude
$$
R_{H_1^n}(\lam) = R_{H_0^n}(\lam)^{1/2} (1+C^n)^{-1} R_{H_0^n}(\lam)^{1/2}
$$
for $\lam$ below the spectrum of $H_0$, where
$$
C^n = (A_0^n R_{H_0^n}(\lam)^{1/2})^* \frac{p_1 - p_0}{p_0} (A_0^n R_{H_0^n}(\lam)^{1/2})
+ R_{H_0^n}(\lam)^{1/2} \frac{q_1 - q_0}{r} R_{H_0^n}(\lam)^{1/2}
$$
and $A_0^n$ denotes the restriction of $A_0$ to $(b_n,b)$ with a Dirichlet boundary
condition at $b_n$.

By assumption, $\frac{p_1 - p_0}{p_0}$ respectively $\frac{q_1 - q_0}{r}$ and thus $\|C^n\|$
can be made arbitrarily small. Hence $(1 +C^n)^{-1}\to 1$ and the first claim follows.

Now, we come to the proof of the relatively nonoscillation part.
Our condition on $p_1/p_0$ imply that
$$
p_0(x)(1 - \eps_-(y)) \leq p_1(x) \leq p_0(x) (1 + \eps_+(y)),
\quad x \geq y,
$$
where
$$
\eps_\pm(y) = \pm\sup_{x \ge y}(\pm (p_1(x)/p_0(x)-1)) \to 0,
\quad y\to b.
$$
Now it follows, from our comparison theorem, that solutions $u_\pm$ of
$\tau_\pm u_\pm =0$ on $(y,b)$, where
$$
\tau_\pm =
\frac 1 r\left(-\frac d{dx} (1 + \eps_\pm(y)) p_0 \frac d{dx} +
q_1 - \lam r \right),
$$
satisfy $\#(u_-,u_0) \geq \#(u_1, u_0) \geq  \#(u_+,u_0)$. Since $u_\pm$
also solve $\ti\tau_\pm u_\pm =0$ on $(y,b)$, where
$$
\ti\tau_\pm =
\frac 1 r\left(-\frac d{dx} p_0 \frac d{dx} +
\frac{q_1 - \lam r}{1 + \eps_\pm(y)} \right),
$$
the result follows from our previous lemma since
$$
\frac{r^{-1} q_1 - \lam}{1 + \eps_\pm(y)} - (r^{-1} q_0 - \lam) =
\frac{r^{-1}(q_1 - q_0) - \eps_\pm(y)(r^{-1} q_0 - \lam)}
{1 + \eps_\pm(y)}
\to 0
$$
as $y\to b$.
\end{proof}

\noindent
Our next aim is to prove Theorem~\ref{thm:reg}. The main ingredient will be Pr\"ufer variables
and the formula (\ref{dwr}) for the derivative of the Wronskian.
Let us suppose that $\tau_{0,1}$ are both regular at
$a$ and $b$ with boundary conditions
\be \label{bc}
\begin{array}{l}
\cos(\alpha) f(a) - \sin(\alpha) p_j(a)f'(a) =0\\
\cos(\beta) f(b) - \sin(\beta) p_j(b)f'(b) =0
\end{array},\quad j=0,1.
\ee
Abbreviate $p_\eps =  p_0 + \eps (p_1 - p_0)$.
Note that $p_\eps^{-1}$ is locally integrable, since $p_\eps^{-1} \le \max(p_0^{-1},p_1^{-1})$.
Hence we can choose $\psi_{\pm,\eps}(\lam,x)$ such that
$\psi_{-,\eps}(\lam,a)=\sin(\alpha)$, $p_\eps(a)\psi_{-,\eps}'(\lam,a)=\cos(\alpha)$ respectively
$\psi_{+,\eps}(\lam,b)=\sin(\beta)$, $p_\eps(b)\psi_{+,\eps}'(\lam,b)=\cos(\beta)$. In particular,
we may choose
\be \label{normtha}
\theta_-(\lam,a) =\alpha \in [0,\pi), \quad -\theta_+(\lam,b) = \pi -\beta \in
[0,\pi).
\ee

Next we introduce
\be
\tau_\eps = \tau_0 + \eps (\tau_1 - \tau_0)=
\frac{1}{r}\left(-\frac{d}{dx}p_\eps\frac{d}{dx}+q_\eps\right),
\quad
\begin{array}{l}
q_\eps =  q_0 + \eps (q_1 - q_0)\\
p_\eps =  p_0 + \eps (p_1 - p_0)
\end{array}
\ee
and investigate the dependence with respect to $\eps\in[0,1]$.

If $u_\eps$
solves $\tau_\eps u_\eps = 0$, then the corresponding Pr\"ufer angles satisfy
\be
\dot{\theta}_\eps(x) = -\frac{W_x(u_\eps, \dot{u}_\eps)}{\rho_\eps^2(x)},
\ee
where the dot denotes a derivative with respect to $\eps$.

As in \cite[Lem.~5.1]{kt}, we obtain by integrating (\ref{dwr}) and using
this to evaluate the corresponding difference quotient the following
lemma.

\begin{lemma} \label{prwpsiepsdot}
We have
\be
W_x(\psi_{\eps,\pm}, \dot{\psi}_{\eps,\pm}) = \left\{
\begin{array}{l} \int_x^b (q_0(t)-q_1(t)) \psi_{\eps,+}(t)^2 dt\\
\qquad {} +\int_x^b (p_1^{-1}(t)-p_0^{-1}(t)) p_\eps \psi_{\eps,+}'(t)^2 dt \\\\
-\int_a^x (q_0(t)-q_1(t)) \psi_{\eps,-}(t)^2 dt \\
\qquad {} + \int_a^x (p_1^{-1}(t)-p_0^{-1}(t)) p_\eps \psi_{\eps,-}'(t)^2 dt, \end{array} \right.
\ee
where the dot denotes a derivative with respect to $\eps$,
$\psi_{\eps,\pm}(x)= \psi_{\eps,\pm}(0,x)$, and
$p_\eps = p_0 + \eps (p_1 - p_0)$.
\end{lemma}

\noindent
Since we assumed $a$ and $b$ to be regular, all integrals exist.

Denote the Pr\"ufer angles of $\psi_{\eps,\pm}(x)= \psi_{\eps,\pm}(0,x)$
by $\theta_{\eps,\pm}(x)$.
The last lemma implies for $q_0-q_1\ge 0$, $p_0-p_1\ge 0$, that
\be
\dot{\theta}_{\eps,+}(x) \le 0, \qquad
\dot{\theta}_{\eps,-}(x) \ge 0.
\ee

Now we are ready to investigate the associated operators $H_0$ and $H_1$.
In addition, we will choose the same boundary conditions for $H_\eps$ as
for $H_0$ and $H_1$. The next lemma follows as in \cite[Lem.~5.2]{kt}.

\begin{lemma} \label{lemevheps}
Suppose $q_0-q_1\ge 0$, $p_0-p_1\ge 0$ (resp.\ both $\le 0$). Then the eigenvalues
of $H_\eps$ are analytic functions with respect to
$\eps$ and they are decreasing (resp.\ increasing).
\end{lemma}

\noindent
In particular, this implies that $\dim\Ran P_{(-\infty,\lam)} (H_\eps)$
is continuous from below (resp.\ above) in $\eps$ for every $\lam$.
Now we are ready for the

\begin{proof}[Proof of Theorem~\ref{thm:reg}]
Without restriction it suffices to assume $\lam_0=\lam_1=0$ and to prove
the result only for $\#(\psi_{0,+}, \psi_{\eps,-})$.

We can split $q_0-q_1$, $p_0-p_1$ in the form
$$
\begin{array}{ll}
q_0 - q_1 = q_+ - q_-,\quad& q_+, q_- \geq 0,\\
p_0 - p_1 = p_+ - p_-,\quad& p_+, p_- \geq 0,
\end{array}
$$
and introduce the operator 
$$
\tau_- = \frac{1}{r}\left(-\frac{d}{dx} (p_0-p_-) \frac{d}{dx} + (q_0-q_-)\right).
$$
Now $\tau_-$ is a negative perturbation of $\tau_0$ and
$\tau_1$ is a positive perturbation of $\tau_-$.

Furthermore, define $\tau_\eps$ by
$$
\tau_\eps = \begin{cases}
\tau_0 + 2\eps (\tau_- - \tau_0), & \eps\in[0,1/2]\\
\tau_- + 2(\eps -1/2)(\tau_1 - \tau_-), & \eps\in[1/2,1].
\end{cases}
$$

Let us look at
$$
N(\eps)=\#(\psi_{0,+}, \psi_{\eps,-}) =
\ceil{\Delta_\eps(b)/\pi} - \floor{\Delta_\eps(a)/\pi} -1, \quad
\Delta_\eps(x)=\Delta_{\psi_{0,+}, \psi_{\eps,- }}(x)
$$
and consider $\eps\in [0,1/2]$.
At the left boundary $\Delta_\eps(a)$ remains constant whereas at the right
boundary $\Delta_\eps(b)$ is increasing by Lemma~\ref{prwpsiepsdot}.
Moreover, it hits a multiple of $\pi$ whenever $0\in\sig(H_\eps)$.
So $N(\eps)$ is a piecewise constant function which is continuous from below
and jumps by one whenever $0\in\sig(H_\eps)$. By Lemma~\ref{lemevheps}
the same is true for
$$
P(\eps) = \dim\Ran\, P_{(-\infty, 0)} (H_\eps) - \dim\Ran\, P_{(-\infty, 0]} (H_0)
$$
and since we have $N(0)=P(0)$, we conclude $N(\eps)=P(\eps)$ for all
$\eps\in[0,1/2]$. To see the remaining case $\eps=[1/2,1]$, simply
replace increasing by decreasing and continuous from below by continuous
from above.
\end{proof}

\section{Approximation in trace norm}

Now we begin with the result for singular
operators by proving the case where
$q_1-q_0$ and $p_1 - p_0$ have compact support.

\begin{lemma} \label{lem:compactcase}
Let $H_j$, $j=0,1$, be Sturm--Liouville
operators on $(a,b)$ associated with $\tau_j$,
and suppose that $r^{-1}(q_1-q_0)$ and $p_1-p_0$
have support in a compact interval $[c,d]\subseteq (a,b)$, where
$a<c$ if $a$ is singular and $d<b$ if $b$ is singular. Moreover, suppose $H_0$ and $H_1$
have the same boundary conditions (if any).

Suppose $\lam_0 < \inf\sig_{ess}(H_0)$. Then
\be
\dim\Ran P_{(-\infty,\lam_0)}(H_1) - \dim\Ran P_{(-\infty,\lam_0]}(H_0)
= \#(\psi_{1,\mp} (\lam_0), \psi_{0,\pm} (\lam_0)).
\ee

Suppose $\sig_{ess}(H_0) \cap [\lam_0,\lam_1] = \emptyset$. Then
\begin{align}
\nonumber
\dim\Ran P_{[\lam_0, \lam_1)} (H_1) &- \dim\Ran P_{(\lam_0, \lam_1]} (H_0)\\
&= \#(\psi_{1,\mp} (\lam_1), \psi_{0,\pm} (\lam_1)) - 
\#(\psi_{1,\mp} (\lam_0), \psi_{0,\pm} (\lam_0)).
\end{align}
\end{lemma}

\begin{proof}
Define $H_\eps= \eps H_1 + (1-\eps) H_0$ as usual and observe that
$\psi_{\eps,-}(z,x)=\psi_{0,-}(z,x)$ for $x\le c$ respectively
$\psi_{\eps,+}(z,x)=\psi_{0,+}(z,x)$ for $x\ge d$. Furthermore,
$\psi_{\eps,\pm}(z,x)$ is analytic with respect to $\eps$ and
$\lam\in\sig_p(H_\eps)$ if and only if $W_d(\psi_{0,+}(\lam),\psi_{\eps,-}(\lam))=0$.
Now the proof can be done as in the regular case. 
\end{proof}

\begin{lemma}\label{lem:compactxi}
Suppose $H_0$, $H_1$ satisfy the same assumptions as in the previous lemma
and that there is a constant $C_1 > 1$ such that $C_1^{-1} \leq p_1(x) p_0(x)^{-1} \leq C_1$
for all $x\in(a,b)$.
Furthermore, set $H_\eps= \eps H_1 + (1-\eps) H_0$. Then
\be
\| \sqrt{r^{-1}\abs{q_0-q_1}} R_{H_\eps}(z) \|_{\mathcal{J}_2} \le C(z), \qquad \eps\in[0,1],
\ee
and
\be
\| \sqrt{\abs{p_1-p_0}} \frac{d}{dx} R_{H_\eps}(z) \|_{\mathcal{J}_2} \le C(z), \qquad \eps\in[0,1].
\ee
In particular, $H_0$ and $H_1$ are resolvent comparable and
\be
\xi(\lambda, H_1, H_0) =  \#(\psi_{1,\mp}(\lam), \psi_{0,\pm}(\lam))
\ee
for every $\lam \in \R \backslash (\sig(H_0) \cup \sig(H_1))$.
Here $\xi(H_1,H_0)$ is assumed to be constructed
such that $\eps \mapsto \xi(H_\eps, H_0)$ is a continuous mapping $[0,1] \rightarrow
L^1((\lam^2+1)^{-1}d\lam)$.
\end{lemma}

\begin{proof}
Denote by
$$
G_\eps(z,x,y) = (H_\eps-z)^{-1}(x,y) = \frac{\psi_{\eps,-}(z,x_<),\psi_{\eps,+}(z,y_>)}{W(\psi_{\eps,-}(z),\psi_{\eps,+}(z))},
$$
where $x_< = \min(x,y)$, $y_>=\max(x,y)$, the Green's function of $H_\eps$. As pointed out
in the proof of the previous lemma, $\psi_{\eps,\pm}(z,x)$ is analytic with respect to $\eps$
and hence a simple estimate shows
$$
\int_a^b \int_a^b |G_\eps(z,x,y)|^2 |r(y)^{-1}(q_1(y)-q_0(y))| r(x)dx\, r(y)dy \le C(z)^2
$$
for $\eps\in[0,1]$, which establishes the first claim.

For the second claim, we need to show that
\begin{align*}
\int_a^b \int_a^b &|p_\eps \partial_x G_\eps(z,x,y)|^2
\left|\frac{p_1(x)-p_0(x)}{p^2_\eps(x)}\right| r(x)dx\, r(y)dy \\
&\leq C(z) \int_c^d \left|\frac{p_1(x)p_0(x)}{p^2_\eps(x)}\right|
\abs{p_0^{-1}(x) -p_1^{-1}(x)} r(x)dx
\end{align*}
is uniformly bounded in $\eps\in[0,1]$. However, this follows
here from the integrand being integrable, since 
$$
0<\frac{p_0}{p_\eps} \leq C_1,\quad
0<\frac{p_1}{p_\eps} \leq C_1.
$$
Moreover, a straightforward calculation (using (\ref{dwr})) and
$$
\psi_{+,\eps} (c) = \psi_{+,\eps'} (c)
W_c(\psi_{+,\eps},\psi_{-,\eps'}) - \psi_{-,\eps'} (c)
W_c(\psi_{+,\eps},\psi_{+,\eps'})
$$
shows
\begin{align*}
G_{\eps'}(z,x,y) = & G_\eps(z,x,y) \\
&+ (\eps-\eps') \int_a^b G_{\eps'}(z,x,t)
r^{-1}(t)(q_1(t)-q_0(t)) G_\eps(z,t,y)  r(t) dt \\
&+ (\eps-\eps') \int_a^b \frac{\partial G_{\eps'}(z,x,t) }{\partial t}
r^{-1}(t)(p_1(t)-p_0(t))
\frac{\partial G_\eps(z,t,y)}{\partial t}  r(t) dt.
\end{align*}
Hence $R_{H_{\eps'}}(z)-R_{H_\eps}(z)$ can be written as the sum
of two products of two Hilbert--Schmidt operators, whose norm can be
estimated by the first claims:
\be
\| R_{H_{\eps'}}(z) - R_{H_\eps}(z) \|_{\mathcal{J}_1} \le |\eps'-\eps| C(z)^2.
\ee
Thus $\eps \mapsto \xi(H_\eps, H_0)$ is continuous. The rest follows from
(\ref{dimxi}).
\end{proof}

\noindent
Before proving Theorem~\ref{thmsing}, we still need to transform
Hypothesis~\ref{hyp:h0h1} in a form such that we can apply our operator
theoretic results from the appendix. The next lemma will do the job.

\begin{lemma}\label{lem:condimp}
Assume Hypothesis~\ref{hyp:h0h1}, and introduce 
\be
\mathfrak{Q} = \{ f \in L^2((a,b),rdx) \,|\, f \in AC_{loc}(a,b), \: \sqrt{p_0} f' \in L^2(a,b)\}.
\ee
Furthermore, introduce the following operators on $\mathfrak{Q}$ with
$N = \lceil \sup_x (p_1(x)p_0(x)^{-1}-1)_+\rceil+1$
\be
A_j = \frac{1}{N^{1/2}} (p_0 - p_1)_+^{1/2} \frac{d}{dx},\quad j = 1,\dots N,
\ee
\be
A_{N+1} = \abs{q_0 - q_1}^{1/2},\quad A_{N+2} = (p_0 - p_1)_-^{1/2} \frac{d}{dx}
\ee
\be
S_1,\dots,S_N=1,\quad S_{N+1}=\sgn(q_0-q_1),\quad S_{N+2} =-1
\ee
Then Hypothesis~\ref{hyp:h0v} is satisfied with these operators and $H_0$, $H_1$ are
self-adjoint extensions of $\tau_0$, $\tau_1$, respectively.
\end{lemma}

\begin{proof}
By Lemma~\ref{lem:infiform}, it is sufficient to check the form
bounds with respect to the form of $\tau_0$ with $q_0 = 0$, since
we have by \cite[Lem.~4.1]{kt}, that $q_0$, $q_1$ will be
infinitesimally form bounded.

To see the claims on the other operators, note that $p_0^{1/2} f'\in L^2$
implies $\abs{p_1 - p_0}^{1/2} f' \in L^2$, since $p_0^{-1/2}\abs{p_1 - p_0}^{1/2}$
is essentially bounded by assumption. We are left with computing the form bounds,
but again ($1\leq j\leq N$, $u\in\mathfrak{Q}$)
$$
\norm{A_j u}^2 =  \frac{1}{N}\norm{p_0^{-1/2}(p_0 - p_1)_+^{1/2} p_0^{1/2} u'}^2 \leq \frac{\sup_x (p_0(x)^{-1}p_1(x)-1)_+}{N} \spr{u}{A_0^* A_0 u}
$$
which shows that the form bound with respect to $A_0^* A_0$ is less then one. By
Lemma~\ref{lem:infiform} the same is true with respect to $H_0$.

Boundedness from below follows by noting, that the quadratic
forms are bounded from below, by the bounds on $q_0$ (resp. $q_1$).
\end{proof}

\noindent
Now we come to the

\begin{proof}[Proof of Theorem~\ref{thmsing}]
We first assume that we have compact support near
one endpoint, say $a$.
Define by $K_\eps$ the multiplication operator by
$\chi_{(a,b_\eps]}$ with $b_\eps \uparrow b$.
Then $K_\eps$ satisfies the assumptions of
Lemma~\ref{lem:resolvconv}.
The last lemma guarantees that Hypothesis~\ref{hyp:h0h1}
implies Hypothesis~\ref{hyp:h0v}, so we can apply Lemma~\ref{lem:Heps} by
Lemma~\ref{lem:resolvconv}.

Denote by 
$$
\tau_\eps = \frac{1}{r}\left(
-\frac{d}{dx} p_\eps \frac{d}{dx}+ q_\eps\right),\quad
\begin{array}{l}
p_\eps = p_0 + \chi_{(a,b_\eps]} (p_1 - p_0) \\
q_\eps = q_0 + \chi_{(a,b_\eps]} (q_1 - q_0)
\end{array}
$$
and by $\psi_{\eps,-}$ the corresponding solutions
satisfying the boundary condition at $a$.
By Lemma~\ref{lem:resolvconv} we have that
$\xi(H_\eps, H_0)$ is constant and equal to $\xi(H_1, H_0)$
once $\eps$ is greater then some $\eps_0$.

Now let us turn to the Wronskians. We first prove the $\#(\psi_{1,-}(\lam), \psi_{0,+}(\lam))$ case.
By Lemma~\ref{lem:compactxi} we know
$$
\xi(\lambda, H_\eps, H_0) = \#(\psi_{\eps,-}(\lam),\psi_{0,+}(\lam)
$$
for every $\eps<1$. Concerning the right-hand side observe that
$$
W_x(\psi_{\eps,-}(\lam), \psi_{0,+}(\lam)) = W_x(\psi_{1,-}(\lam), \psi_{0,+}(\lam))
$$
for $x \le b_\eps$ and that $W_x(\psi_{\eps,-}(\lam), \psi_{0,+}(\lam))$ is constant for
$x \ge b_\eps$. This implies that for $\eps\ge \eps_0$ we have
\begin{align*}
\xi(\lam, H_1, H_0) &= \xi(\lam, H_\eps, H_0) = \#(\psi_{\eps,-}(\lam), \psi_{0,+}(\lam))\\
& = \#_{(a,b_\eps)}(\psi_{\eps,-}(\lam), \psi_{0,+}(\lam)) = \#_{(a,b_\eps)}(\psi_{1,-}(\lam), \psi_{0,+}(\lam)).
\end{align*}
In particular, the last item $\#_{(a,b_\eps)}(\psi_{1,-}(\lam), \psi_{0,+}(\lam))$ is eventually constant
and thus has a limit which, by definition, is $\#(\psi_{1,-}(\lam), \psi_{0,+}(\lam))$.

For the corresponding $\#(\psi_{1,+}(\lam), \psi_{0,-}(\lam))$ case one simply exchanges the
roles of $H_0$ and $H_1$.

Hence the result holds if the perturbation has compact support near one endpoint. Now one repeats the
argument to remove the compact support assumption near the other endpoint as well.
\end{proof}

\section{Appendix: The Spectral Shift Function}
\label{sec:xi}

In this appendix we collect some facts 
on Krein's spectral shift function which are of
relevance to us. Most results are taken from
\cite{yafams} (see also \cite{sro} for an easy introduction).
The first part closely follows the appendix of \cite{kt}.

Two operators $H_0$ and $H_1$ are called resolvent comparable, if
\be\label{eq:resolvcond}
R_{H_1}(z) - R_{H_0} (z)
\ee
is trace class for one $z\in\rho(H_1)\cap\rho(H_0)$. By the first resolvent identity \eqref{eq:resolvcond}
then holds for all $z\in\rho(H_1)\cap\rho(H_0)$.

\begin{theorem}[Krein \cite{krein}]\label{thm:shifttheorem}
Let $H_1$ and $H_0$ be two resolvent comparable self-adjoint operators, 
then there exists a function
\be
\xi(\lam, H_1, H_0) \in L^1(\R, (\lambda^2 + 1)^{-1}d\lam)
\ee
such that 
\be\label{eq:traceformula}
\tr(f(H_1) - f(H_0)) = \int_{-\infty}^\infty \xi(\lam, H_1, H_0) f'(\lambda) d\lam
\ee
for every smooth function $f$ with compact support.
\end{theorem}
Note: Equation \eqref{eq:traceformula} holds in fact for a much larger class of functions $f$.
See \cite[Thm.~9.7.1]{yafams} for this and a proof of the last theorem.

The function $\xi(\lam) = \xi(\lam, H_1, H_0)$ is called Krein's spectral shift function and is
unique up to a constant. Moreover, $\xi(\lam)$ is clearly constant on every interval
$(\lam_0, \lam_1) \subset \rho(H_0)\cap\rho(H_1)$.
Hence, if $\dim\Ran P_{(\lam_0, \lam_1)}(H_j) < \infty$, $j = 0,1$,
then $\xi(\lam)$ is a step function and
\be \label{dimxi}
\dim \Ran P_{(\lam_0, \lam_1)} (H_1) - \dim \Ran P_{(\lam_0, \lam_1)} (H_0) =
\lim_{\eps\downarrow0} \Big( \xi(\lam_1-\eps) - \xi(\lam_0+\eps) \Big).
\ee
This formula clearly explains the name spectral shift function.

Before investigating further the properties of the SSF,
we will recall a few things about trace ideals
(see for example \cite{str}). First,
for $1\leq p <\infty$ denote by $\mathcal{J}^p$ the
Schatten $p$-class, and by $\norm{.}_{\mathcal{J}^p}$ its
norm. We will use $\norm{.}$ for the usual operator
norm. Using $\norm{A}_{\mathcal{J}^p} = \infty$ if
$A \notin \mathcal{J}^p$, we have the following
inequalities for all operators:
$$
\norm{A B}_{\mathcal{J}^p} \leq
\norm{A}\norm{B}_{\mathcal{J}^p},\quad
\norm{A B}_{\mathcal{J}^1} \leq \norm{A}_{\mathcal{J}^2}
\norm{B}_{\mathcal{J}^2}.
$$
Furthermore, we will use the notation of
$\mathcal{J}^p$-converges to denote convergence in the
respective norm.
The following result from \cite[Thm.~IV.11.3]{tracedet} will be needed.

\begin{lemma}\label{lem:contrace}
Let $p > 0$, $A_n \xrightarrow{\mathcal{J}^p} A$, $T_n \xrightarrow{s} T$, $S_n \xrightarrow{s} S$ sequences
of strongly convergent bounded linear operators, then:
\be
\|T_n A_n S_n^\ast - T A S^\ast\|_{\mathcal{J}^p} \rightarrow 0.
\ee
Here $\|.\|_{\mathcal{J}^p}$ are the norms of the Schatten $p$-classes $\mathcal{J}^p$.
\end{lemma}

\noindent
We will also need the following continuity result for $\xi$.
It will also allow us to fix the unknown constant. The second
part is \cite[Lem.~7.3]{kt}, the first from \cite{yafams}.

\begin{lemma} \label{lem:Heps}
Suppose $H_\eps$, $\eps\in[0,1]$, is a family of self-adjoint operators, which is continuous in
the metric
\be
\rho(A,B) = \|R_A (z_0) - R_B(z_0)\|_{\mathcal{J}^1}
\ee
for some fixed $z_0\in \C\backslash\R$ and abbreviate $\xi_\eps = \xi(H_\eps, H_0)$.
Then there exists a unique choice of $\xi_\eps$ such that $\eps \mapsto \xi_\eps$ is continuous
$[0,1] \rightarrow L^1(\R,(\lam^2+1)^{-1}d\lam)$ with $\xi_0 = 0$.

If $H_\eps \ge \lam_0$ is bounded from below, we can also allow $z_0 \in (-\infty,\lam_0)$.

For $\lam\in\rho(H_1)\cap\R$, we have that there is an $\eps_0$ such that
$\xi_\eps(\lam)=\xi_1(\lam)$ for $\eps>\eps_0$.
\end{lemma}

\begin{proof}
We just need to proof the third part.
For $\eps$ close to $1$ a whole neighborhood of $\lam$
is in $\rho(H_\eps)\cap\R$, since the resolvent sets converge.
Furthermore, we know from this that the $\xi_\eps$ is integer
valued near $1$ in a neighborhood of $\lam$.
Now the claim follows from the convergence of $\xi_\eps\to\xi_1$
in $L^1(\R,(\lam^2+1)^{-1}d\lam)$.
\end{proof}

\noindent
Our final aim is to find some conditions which
allow us to verify the assumptions of this
lemma.  To do this, we derive some properties of relatively bounded operators
multiplied by strongly continuous families of operators.

\begin{hypothesis} \label{hyp:h0v}
Suppose $H_0$ is self-adjoint and bounded from below.
Let $A_j$, $j=1,\dots,n$, be closed operators and $S_j$, $j=1,\dots,n$
be bounded operators with $\norm{S_j}\leq 1$.
Furthermore, suppose that these satisfy for $j=1,\dots,n-1$ that
\begin{enumerate}
\item[(i)] $A_j^\ast A_j$ is relatively form bounded with respect
to $H_0$ with relative form bound less than one and $S_j$ is positive, or
\item[(i')] $A_j^\ast A_j$ is infinitesimally form bounded with respect to $H_0$.
\end{enumerate}
Suppose for $j=n$, that
\begin{enumerate}
\item[(ii)] $A_j^\ast A_j$ is relatively form bounded with respect
to $H_0$ with relative form bound less than one.
\end{enumerate}
\end{hypothesis}

\noindent
Note that condition (i) implies that $A_j^\ast S_j A_j$ is a positive operator.

We recall that $A^\ast A$ being form bounded with respect to $H_0$
means that we have $\fdom(A^\ast A) \supseteq \fdom(H_0)$ and
\be \label{defrfb}
\spr{\psi}{A^\ast A \psi} \leq a \spr{\psi}{H_0 \psi} + b \|\psi\|^2,
\quad \forall \psi \in \fdom(H_0).
\ee
for some $0 \leq a < 1$, $0 \leq b$. The form bound is the infinimum over
all $a$ such that (\ref{defrfb}) holds.

The next lemma is modified from \cite[Lem.~7.5]{kt},
to be able to deal with differential operators and sums
of operators.

\begin{lemma}\label{lem:onetermHeps}
Let $\eps \ni[0,1] \to K_\eps$ be a strongly
continuous family of self-adjoint bounded operators
which satisfy $0= K_0\le K_\eps \le K_1 = 1$.

Let 
\begin{enumerate}
\item
$\eps\mapsto H_\eps$ satisfy the assumptions
of Lemma~\ref{lem:Heps},
\item
$S$ be a bounded operator with $\norm{S}\leq 1$, and
\item
$A$ be a closed operator such that
$A^\ast A$ is relatively  bounded with respect to
$H_\eps$ with uniform in $\eps$ bound less then one, 
and $A R_{H_\eps}(z)\in\mathcal{J}^2$
for one $z\in\C\backslash\R$.
\end{enumerate}
Then $\ti{H}_\eps = H_\eps + A^\ast K_\eps S A$
also satisfies the assumptions of Lemma~\ref{lem:Heps}.
Furthermore, for form bounded $B$ with $B R_{H_\eps}(z)\in\mathcal{J}^2$,
we have $B R_{\ti{H}_\eps}(z)\in\mathcal{J}^2$ for all $\eps\in[0,1]$.
\end{lemma}

\begin{proof}
We will abbreviate $V_\eps = A^\ast K_\eps S A$,
$\tilde{H}_\eps = H_\eps + V_\eps$, $R_\eps(z)= R_{H_\eps}(z)$,
and $\ti{R}_\eps(z)= R_{\ti{H}_\eps}(z)$.
By the KLMN theorem (\cite[Thm.~X.17]{reedsimon2}), $\ti{H}_\eps$
is self-adjoint since
$$
|\spr{\psi}{V_\eps \psi}| \leq
|\spr{A\psi}{K_\eps SA\psi}|
\le \spr{\psi}{A^\ast A \psi}, \quad
\psi\in\fdom(V_\eps).
$$
Moreover, using (\ref{defrfb}) we obtain
$$
\| A R_\eps(-\lam)^{1/2} \|^2 \le a, \qquad \frac{b}{a}< \lam.
$$
For $\lam>\frac{b}{a}$ we have by Lemma~\ref{lem:secresoform}
\begin{align*}
\ti{R}_\eps(-\lam) &= R_\eps(-\lam)^{1/2} (1+C_\eps)^{-1} R_\eps(-\lam)^{1/2}, \\
C_\eps &= (A R_\eps(-\lam)^{1/2})^\ast (K_\eps S A R_\eps(-\lam)^{1/2}).
\end{align*}
Hence, a straightforward calculation shows
\begin{align}\label{eq:formsecresol}
\ti{R}_\eps(-\lam) &= R_\eps(-\lam) - (A R_\eps(-\lam))^*
(1 +\ti{C}_\eps)^{-1} (K_\eps S A R_\eps(-\lam)),\\
\nn \ti{C}_\eps &= (K_\eps S A R_\eps(-\lam)^{1/2})
(A R_\eps(-\lam)^{1/2})^\ast.
\end{align}
By $\|\ti{C}_\eps\| \le a < 1$, we have that
$(1 + \ti{C}_\eps)^{-1}$ exists. Furthermore,
note that (\ref{eq:formsecresol}) implies, that
$B \ti{R}_\eps(-\lam) \in \mathcal{J}^2$, since:
$$
B \ti{R}_\eps(-\lam) = B R_\eps(-\lam) - 
B R_\eps(-\lam)^{1/2} (A R_\eps(-\lam)^{1/2})^*
(1 +\ti{C}_\eps)^{-1} (K_\eps S A R_\eps(-\lam)),
$$
and $A R_\eps(-\lam)\in\mathcal{J}^2$.
Now, look at
\begin{align*}
D_{\eps,\eps'}\psi &= (- \ti{C}_\eps(1 + \ti{C}_\eps)^{-1} -
\ti{C}_{\eps'} (1 + \ti{C}_{\eps'})^{-1})\psi \\
\nn &= (C_{\eps'} - C_\eps)(1 + C_\eps)^{-1} \psi -
C_{\eps'} D_{\eps,\eps'} \psi,
\end{align*}
where
$$
D_{\eps,\eps'} = (1 + \ti{C}_\eps)^{-1} - (1 + \ti{C}_{\eps'})^{-1}.
$$
Taking norms we obtain
$$
\|D_{\eps,\eps'}\psi\| = \frac{1}{1 - a} \|(C_{\eps'} - C_\eps)(1 + C_\eps)^{-1} \psi\|,
$$
where the last term converges to $0$ as $\eps'\to\eps$.
This implies, that $(1 + \ti{C}_\eps)^{-1}$ is strongly
continuous. Now, we obtain from (\ref{eq:formsecresol}) for the difference of resolvents
\begin{align*}
\ti{R}_\eps(-\lam) - \ti{R}_{\eps'}(-\lam)
= (A R_\eps(-\lam))^* ((1 +\ti{C}_\eps)^{-1} K_\eps - 
(1 +\ti{C}_{\eps'})^{-1} K_{\eps'}))
 (S A R_\eps(-\lam))
\end{align*}
$\mathcal{J}^1$-converges to $0$ as $\eps\to\eps'$ by
Lemma~\ref{lem:contrace} and by
$A R_\eps(-\lam)\in\mathcal{J}^2$.  This way we also obtain that
$\ti{H}_\eps$ and $\ti{H}_{\eps'}$ are indeed resolvent comparable.
\end{proof}

\noindent
We also recall the following well-known fact on quadratic forms:

\begin{lemma}\label{lem:infiform}
Let $v,s,t$ be quadratic forms, such that $s$ is positive
and symmetric, and $v$ is infinitesimal form bounded
with respect to $s$, and $t$ is form bounded with bound less
then $1$ with respect to $s$. Then $t$ is also form
bounded with bound less then $1$ with respect to $s + v$.
\end{lemma}

\begin{proof}
Using $|v(\psi)| \le \eps s(\psi) + C\|\psi\|^2$ for arbitrary small $\eps > 0$,
a direct calculation shows
$$
s(\psi) \leq \frac 1{1 - \eps} \abs{s(\psi) + v(\psi)}
+ \frac{C(\eps)}{1 - \eps}\norm{\psi}^2.
$$
Denoting by $a$ the $s$ bound of $t$, it follows that
$t$ is $s + v$ bounded with bound less then $a/(1 - \eps)$,
implying that the bound is again less then one.
\end{proof}

\begin{lemma}\label{lem:resolvconv}
Let $\eps \ni[0,1] \to K_\eps$ be a strongly
continuous family of self-adjoint bounded operators
which satisfy $0= K_0\le K_\eps \le K_1 = 1$. 

Assume Hypothesis~\ref{hyp:h0v}. Then
\be
H_\eps = H_0 + \sum_{j=1}^n A_j ^\ast K_\eps S_j A_j
\ee
are self-adjoint operators such that the assumptions of Lemma~\ref{lem:Heps}
hold.
\end{lemma}

\begin{proof}
Introduce $H_\eps^m = H_0 + V_\eps^m$, $m = 0, \dots, n$, where
$V_\eps^m = \sum_{j=1}^m A_j^\ast K_\eps S_j A_j$.
Since all but the last perturbations are either positive
or infinitesimal (in which case on has to use Lemma~\ref{lem:infiform}),
we can assume that $A_m^\ast K_\eps S_m A_m$
is relatively form bounded with uniform bound less then one with
respect to $H_\eps^l$ with $l<m$.

Now, the result follows by applying the previous lemma
with $H = H^{m-1}$, $\ti{H} = H^m$, $A = A_m$, $S = S_m$
and $B = B_l$, $l = m+1,\dots,n$ and letting $m$ going up
from $1$ to $n$.
\end{proof}

\end{document}